\documentclass[10pt]{article}

\usepackage{amscd,amsmath, amssymb, fancyhdr, graphicx, url}

\usepackage{dutchcal}

\numberwithin{equation}{section}

\newcommand{\version}{version 3.0,\ \ 11.07.2017}

\makeatletter
\def\x@arrow{\DOTSB\Relbar}
\def\xlongrightarrowfill@{\arrowfill@\relbar\relbar\longrightarrow}
\newcommand{\xlongrightarrow}[2][]{%
        \ext@arrow 0099\xlongrightarrowfill@{#1}{#2}}
\makeatother

\def\eqref#1{(\ref{#1})}

\newcommand{\arrow}{{\:\longrightarrow\:}}
\newcommand{\Z}{{\Bbb Z}}
\def\C{{\Bbb C}}

\newcommand{\Char}{{\operatorname{\sf char}}}

\newcommand{\R}{{\Bbb R}}
\newcommand{\Q}{{\Bbb Q}}

\def\1{\sqrt{-1}\:}

\newcommand{\cntrct}                
{\hspace{2pt}\raisebox{1pt}{\text{$\lrcorner$}}\hspace{2pt}}

\renewcommand{\tilde}{\widetilde}
\renewcommand{\bar}{\overline}
\renewcommand{\phi}{\varphi}
\renewcommand{\epsilon}{\varepsilon}
\renewcommand{\geq}{\geqslant}
\renewcommand{\leq}{\leqslant}

\newcommand{\Perspace}{\operatorname{{\Bbb P}\sf er}}

\newcommand{\prim}{\text{\sf prim}}

\newcommand{\End}{\operatorname{End}}
\newcommand{\MT}{\operatorname{MT}}

\newcommand{\St}{\operatorname{St}}

\newcommand{\Gal}{\operatorname{{\mathcal{G\!\!\:a\!\!\:l}}}}
\newcommand{\Cl}{\operatorname{Cl}}

\newcommand{\Sym}{\operatorname{Sym}}
\newcommand{\Pic}{\operatorname{Pic}}

\newcommand{\Aut}{\operatorname{Aut}}

\newcommand{\proof}{{\bf Proof:\ }}

\newcounter{Mycounter}[section]
\newcounter{lemma}[section]
\newcounter{claim}[section]
\renewcommand{\theclaim}{{Claim \thesection.\arabic{claim}}}
\newcommand{\claim}{%
    \setcounter{claim}{\value{Mycounter}}
    \refstepcounter{claim}
    \stepcounter{Mycounter}
    {\noindent \bf \theclaim:\ }}

\newcounter{sublemma}[section]
\newcounter{corollary}[section]
\renewcommand{\thecorollary}{{Corollary \thesection.\arabic{corollary}}}
\newcommand{\corollary}{%
    \setcounter{corollary}{\value{Mycounter}}
    \refstepcounter{corollary}
    \stepcounter{Mycounter}
    {\noindent \bf \thecorollary:\ }}

\newcounter{theorem}[section]
\renewcommand{\thetheorem}{{Theorem \thesection.\arabic{theorem}}}
\newcommand{\theorem}{%
    \setcounter{theorem}{\value{Mycounter}}
    \refstepcounter{theorem}
    \stepcounter{Mycounter}
    {\noindent \bf \thetheorem:\ }}

\newcounter{conjecture}[section]
\newcounter{proposition}[section]
\renewcommand{\theproposition}
      {{Proposition \thesection.\arabic{proposition}}}
\newcommand{\proposition}{%
    \setcounter{proposition}{\value{Mycounter}}
    \refstepcounter{proposition}
    \stepcounter{Mycounter}
    {\noindent \bf \theproposition:\ }}

\newcounter{definition}[section]
\renewcommand{\thedefinition}
      {{Definition~\thesection.\arabic{definition}}}
\newcommand{\definition}{%
    \setcounter{definition}{\value{Mycounter}}
    \refstepcounter{definition}
    \stepcounter{Mycounter}
    {\noindent \bf \thedefinition:\ }}

\newcounter{example}[section]
\renewcommand{\theexample}{{Example \thesection.\arabic{example}}}
\newcommand{\example}{%
    \setcounter{example}{\value{Mycounter}}
    \refstepcounter{example}
    \stepcounter{Mycounter}
    {\noindent \bf \theexample:\ }}

\newcounter{remark}[section]
\renewcommand{\theremark}{{Remark \thesection.\arabic{remark}}}
\newcommand{\remark}{%
    \setcounter{remark}{\value{Mycounter}}
    \refstepcounter{remark}
    \stepcounter{Mycounter}
    {\noindent \bf \theremark:\ }}

\newcounter{problem}[section]
\newcounter{question}[section]
\makeatletter

\@addtoreset{equation}{section} \@addtoreset{footnote}{section}
\makeatother

\def\blacksquare{\hbox{\vrule width 5pt height 5pt depth 0pt}}
\def\endproof{\blacksquare}

\begin{document}

\begin{center}
{\LARGE\bf
Transcendental Hodge algebra\\[3mm]
}

Misha
Verbitsky\footnote{Misha Verbitsky is partially supported by the 
Russian Academic Excellence Project '5-100'.


 }

\end{center}

{\small \hspace{0.15\linewidth}
\begin{minipage}[t]{0.7\linewidth}
{\bf Abstract} \\
The transcendental Hodge lattice of a projective manifold $M$
is the smallest Hodge substructure in $p$-th cohomology which
contains all holomorphic $p$-forms. We prove that the direct
sum of all transcendental Hodge lattices has a natural
algebraic structure, and compute this algebra explicitly
for a hyperk\"ahler manifold. As an application, we obtain
a theorem about dimension of a compact torus $T$
admitting a holomorphic symplectic embedding to a hyperk\"ahler
manifold $M$. If $M$ is generic in a $d$-dimensional 
family of deformations, then $\dim T\geq 2^{[(d+1)/2]}$. 
\end{minipage}
}

\tableofcontents


\section{Introduction}


Transcendental Hodge lattice is the smallest substructure
in a Hodge structure of a projective manifold containing
the cohomology classes of all holomorphic $p$-forms.

Transcendental Hodge lattices were (as Yu. Zarhin notices)
a somewhat neglected subject in Hodge theory. 
This is somewhat surprising, because, unlike the
Hodge structure on cohomology, the transcendental
Hodge lattice is a birational invariant.

The observation made in this paper was probably made
before, but it is still very useful: the direct sum
of all transcendental Hodge lattices is naturally
an algebra (Subsection \ref{_tra_alge_Subsection_}).
Using Zarhin's classification of transcendental 
Hodge structures of K3 type, we compute this
algebra for all hyperk\"ahler manifolds
(\ref{_transcende_Theorem_}; see Subsection \ref{_hk_Subsection_}
for the definition and basic properties of hyperk\"ahler manifolds).
This computation gives a way to several generalization 
of results which were previously known for general
(non-algebraic) hyperk\"ahler manifolds, proving non-existence
of low-dimensional holomorphic symplectic tori in projective 
hyperk\"ahler manifolds (\ref{_k_symple_in_families_Corollary_}). 

We expect more results to be obtained in the same direction,
because the transcendental Hodge algebra has a promise to become
a very powerful tool in the study of projective holomorphically
symplectic varieties.

\subsection{Mumford-Tate group and Hodge group}

Mumford defined the Mumford-Tate group
in \cite{_Mumford:families_}, and called in ``Hodge
group'', and in literature these terms are sometimes
considered equivalent. However, in Zarhin's papers 
\cite{_Zarhin:Hodge_K3_,_Zarhin:Hodge_K3_LNM_,_Zarhin:small_number_}
as in \cite{_Peters_Steenbrink:MHS_}, these notions are distinct.
We shall call it ``special Mumford-Take group''.

Given a Hodge structure $V=\bigoplus V^{p,q}$,
we consider action of $\C^*$ on $V$ where $z$ acts
on $V^{p,q}$ as multiplication by $z^{p-q}$.
Then ``Mumford-Tate group'' is the smallest rational
algebraic group containing image of $\C^*$ and
``special Mumford-Tate group'', or
``Hodge group'', is the smallest rational
algebraic group containing image of $U(1)\subset \C^*$

We shall not pay much attention to this difference.


\subsection{Hyperk\"ahler manifolds: an introduction}
\label{_hk_Subsection_}

Here we list some basic facts and properties of hyperk\"ahler
manifolds; for more details, proofs and history, please see 
\cite{_Besse:Einst_Manifo_} and \cite{_Beauville_}.

\hfill

\definition
A {\bf hyperk\"ahler structure} on a manifold $M$
is a Riemannian structure $g$ and a triple of complex
structures $I,J,K$, satisfying quaternionic relations
$I\circ J = - J \circ I =K$, such that $g$ is K\"ahler
for $I,J,K$.

\hfill

\remark Let $\omega_I, \omega_J, \omega_K$
be the K\"ahler symplectic forms associated with $I,J,K$:
\[
\omega_I(\cdot, \cdot)=g(\cdot, I\cdot), \ \ \omega_J(\cdot, \cdot)=g(\cdot, J\cdot),\ \ 
\omega_K(\cdot, \cdot)=g(\cdot, K\cdot).
\]
A hyperk\"ahler manifold is holomorphically
symplectic: $\omega_J+\1 \omega_K$ is a holomorphic
symplectic form on $(M,I)$. Converse is also true:

\hfill

\theorem (Calabi-Yau; see \cite{_Yau:Calabi-Yau_}, \cite{_Besse:Einst_Manifo_}) 
A compact, K\"ahler, holomorphically symplectic manifold
admits a unique hyperk\"ahler metric in any K\"ahler class.

\hfill

\definition A compact hyperk\"ahler manifold $M$ is called
{\bf maximal holonomy manifold}, or {\bf simple},
or {\bf IHS} if $\pi_1(M)=0$, $H^{2,0}(M)=\C$.

\hfill

\theorem
(Bogomolov Decomposition Theorem; \cite{_Bogomolov:decompo_}.) \\
Any hyperk\"ahler manifold admits a finite covering
which is a product of a torus and several 
maximal holonomy (simple) hyperk\"ahler manifolds.
\endproof

\hfill

\remark
From now on, all holomorphic symplectic
manifolds are tacitly assumed to be of K\"ahler type,
``holomorphic symplectic'' is used interchangeably
with ``hyperk\"ahler'', and all hyperk\"ahler manifolds
are assumed to be of maximal holonomy.

\subsection{Trianalytic and holomorphic symplectic subvarieties}

The starting point of the study of holomorphically
symplectic subvarieties in hyperk\"ahler manifolds
was the paper \cite{_Verbitsky:Symplectic_I_}, 
where it was shown that any complex subvariety of a general
(non-algebraic) hyperk\"ahler manifold is holomorphically
symplectic outside of its singularities. In 
\cite{_Verbitsky:Symplectic_II_} this result 
was improved: it was shown that for a general
complex structure induced by a unit quaternion
$L=aI+bJ +cK$, $L^2=-1$ on a hyperk\"ahler manifold
$(M,I,J,K)$, any complex subvariety of $(M,L)$
is in fact {\bf trianalyric}, that is, 
complex analytic with respect to $I, J, K$.
The notion of a trianalytic subvariety, introduced in this paper, 
had many uses further on. In \cite{_Verbitsky:hypercomple_},
the notion of an abstract ``hypercomplex manifold'' was developed,
following Deligne and Simpson
(\cite{_Deligne:defi_}, 
\cite{_Simpson:hyperka-defi_}). The examples of hypercomplex
varieties include all trianalytic subvarieties
of hyperk\"ahler manifold. It was shown that
hypercomplex varieties admit a natural 
hypercomplex desingularization. Applied to trianalytic
subvarieties of $M$, this gives a holomorphically
symplectic desingularization immersed to $M$
holomorphically symplectically.

In \cite{_Verbitsky:non-compact_subva_}, 
the results of \cite{_Verbitsky:Symplectic_II_}
were  generalized to non-compact hyperk\"ahler
manifolds, with further generalizations in 
\cite{_Soldatenkov_Verbitsky:subvarieties_}.

Much advance was done towards classifying trianalytic subvarieties
in general deformations of hyperk\"ahler manifolds
(\cite{_Verbitsky:Hilbert_}, \cite{_KV:partial_}, 
\cite{_Kurnosov:tori_}, \cite{_SV:k-symplectic_}). 
However, the subject of more general (that is,
not necessarily trianalytic) holomorphic symplectic
subvarieties was mostly neglected, with the paper
\cite{_V:Wirtinger_} being the only exception (as far as we know).

In \cite{_V:Wirtinger_}, the 
{\bf Wirtinger invariant}$\mu(Z,M)$ of a holomorphically
symplectic (possibly, immersed) subvariety $Z$ in
a hyperk\"ahler manifold $M$ was defined.
Wirtinger number measures how far a holomorphically
symplectic subvariety $Z\subset M$  is from being 
trianalytic. This invariant satisfies the inequality
$\mu(Z,M) \geq 1$, with equality reached if and only if
$Z$ is trianalytic.

It was shown that $\mu(Z,M)$ is monotonous and multiplicative,
in the following sense.
Given a chain of symplectic immersions  $Z_1\arrow Z_2\arrow M$, one has
\[ \mu(Z_1, M) =\mu(Z_1, Z_2) \mu(Z_2, M).\] Therefore, $Z_1$ is trianalytic
in $M$ if and only if it is trianalytic in $Z_2$ and $Z_2$ is
trianalytic in $M$.

The formalism of transcendental Hodge algebra seems to be particularly
suited for the  study of the holomorphically symplectic subvarieties. In this
paper, we generalize the results of \cite{_SV:k-symplectic_} from trianalytic
to more general symplectic subvarieties. In \cite{_SV:k-symplectic_} 
it was shown that a trianalytic complex subtorus $Z$ in a very general deformation
of a hyperk\"ahler manifold $M$ satisfies $\dim Z\leq 2^{\left\lfloor \frac{d+1}2\right\rfloor}$,
where $d=b_2-2$ is the dimension of the universal family of deformations of
$M$. In this paper, the same result is proven for projective $M$
generic in a $d$-dimensional family of deformations (\ref{_k_symple_in_families_Corollary_}).


\section{Hodge structures}


In this section we briefly introduce the Hodge structures;
for more examples, context and applications, see
\cite{_Voisin-Hodge_}  and \cite{_Griffiths:transcendental_}.

\subsection{Hodge structures: the definition}

\definition
Let $V_\R$ be a real vector space.
{\bf A (real) Hodge structure of weight $w$} 
on a vector space $V_\C=V_\R \otimes_\R \C$ 
is a decomposition $V_\C =\bigoplus_{p+q=w} V^{p,q}$, satisfying 
$\overline{V^{p,q}}= V^{q,p}$. It is called {\bf rational
Hodge structure} if one fixes a rational lattice $V_\Q$
such that $V_\R=V_\Q \otimes \R$. A Hodge structure is 
equipped with $U(1)$-action, with $u\in U(1)$
acting as $u^{p-q}$ on $V^{p,q}$. {\bf Morphism}
of Hodge structures is a rational map which is $U(1)$-invariant.
Notice that this $U(1)$-action is complex conjugate
to itself. This means that it is well defined on $V_\R$.

Throughout this paper we often  use
``Hodge structures'' as a shorthand for
``rational Hodge structures''.

\hfill

\definition 
{\bf Polarization}
on a rational Hodge structrure of weight $w$ is 
non-degenerate 2-form $h\in V_\Q^*\otimes V^*_\Q$ 
(symmetric or antisymmetric depending on parity of $w$) which 
satisfies 
\begin{equation}\label{_Hodge_Riemann_Equation_}
 -(\1)^{p-q}h(x, \bar x)>0 
\end{equation} (``Riemann-Hodge relations'')
for each non-zero $x\in V^{p,q}$, and
$h(x, y)=0$ for any $x\in V^{p,q}$,
$y\in V^{p',q'}$, unless $p\neq q', q\neq p'$.
The later condition is equivalent to $U(1)$-invariance
of $h$ on $V_\R$.

\hfill

\remark
Further on, the Hodge structures we consider are
tacitly assumed to be rational and polarized.

\hfill

\example 
Let $(M, \omega)$ be a compact K\"ahler manifold. Then the Hodge decomposition
$H^*(M,\C)=\oplus H^{p,q}(M)$ defines a Hodge structure on $H^*(M,\C)$.
If we restrict ourselves to the primitive cohomology space
\[ H^*_\prim:= \{ \eta\in H^*(M) \ \ (*\eta)\wedge\omega=0\},
\]
and consider $h(x,y):= \int_M x\wedge y \wedge \omega^{\dim_\C M -w}$,
relations \eqref{_Hodge_Riemann_Equation_} become the usual Hodge-Riemann relations.
If, in addition, the cohomology class of $\omega$ is rational
(in this case, by Kodaira theorem, $(M,\omega)$ is projective)
the space $H^*_\prim(M)$ is also rational, and the Hodge decomposition
$H^*_\prim(M,\C)=\oplus H^{p,q}_\prim(M)$ defines a 
polarized, rational Hodge structure.

\subsection{Special Mumford-Tate group}
\label{_MT_Subsection_}

Results of Subsection \ref{_MT_Subsection_} and 
Subsection \ref{_MT_Gal_Subsection_} are rather standard;
please see \cite{_Moonen:notes_} and third part of
\cite{_Zarhin:weights_}. We repeat these results here
for further use.

\hfill

\definition
A {\bf simple object} of an abelian category is an object
which has no proper subobjects.
An abelian category is {\bf semisimple} if any object is
a direct sum of simple objects.

\hfill

\claim
 Category of polarized Hodge structures in semisimple

\hfill

\proof 
Orthogonal complement of a Hodge substructure $V'\subset V$
with respect to $h$ is again a Hodge substructure, and
this complement does not intersect $V'$; both assertions
follow from the Riemann-Hodge relations. \endproof

\hfill

\definition
Let $V$ be a Hodge structure over $\Q$,
and $\rho$ the corresponding $U(1)$-action.
{\bf Special Mumford-Tate group} (Mumford, \cite{_Mumford:families_}; 
Mumford called it ``the Hodge group'') is the
smallest algebraic group over $\Q$ containing
$\rho$.

\hfill

\theorem\label{_Mumford_Tate_via_invariants_Theorem_}
Let $V$ be a rational, polarized Hodge structure, and $\MT(V)$
its special Mumford-Tate group. Consider the tensor algebra
of $V$, $W= T^{\otimes}(V)$ with the Hodge structure
(also polarized) induced from $V$. Let $W_h$ be the 
space of all $\rho$-invariant rational vectors in $W$
(such vectors are called {\bf ``Hodge vectors''}). 
Then $\MT(V)$
coincides with the stabilizer 
\[ \St_{GL(V)}(W_h):= \{ g\in GL(V)\ \ |\ \ \forall w\in W_h, g(w)=w\}.
\]
\proof Follows from the Chevalley's theorem on tensor invariants;
however, to apply it one needs first to show that the
Mumford-Tate group is reductive.
\endproof

\hfill

\corollary\label{_MT_inva_=_Hodge_substra_Corollary_}
Let $V_\Q$ be a rational Hodge structure, and $W\subset V_\C$
a subspace. Then the following are equivalent:
\begin{description}
\item[(i)] $W$ is a Hodge substructure.
\item[(ii)] $W$ is special Mumford-Tate invariant.
\end{description}
\endproof

\subsection{Special Mumford-Tate group and the
  $\Aut(\C/\Q)$-action}
\label{_MT_Gal_Subsection_}

\definition
 Let $V_\Q$ be a $\Q$-vector space, and $V_\C:= V_\Q\otimes_\Q \C$
its complexification, equipped with a natural
$\Aut(\C/\Q)$ action. We call a complex subspace $T\subset V_\C$
{\bf rational} if $T=T_\Q \otimes_\Q \C$, where
$T_Q= V_\Q \cap T$.

\hfill

\remark 
A complex subspace $W\subset V_\C$ is rational
if and only if it is $\Aut(\C/\Q)$-invariant.

\hfill

This implies

\hfill

\claim\label{_Gal_gene_smallest_HS_Claim_}
Consider a subspace $W\subset V_\C$, and
let $\tilde W_\Q\subset V_\Q$ be the smallest subspace of $V_\Q$
such that $\tilde W_\Q\otimes_\Q \C\supset W$. Then $W_\Q\otimes_\Q \C$
is generated by $\sigma(W)$, for all $\sigma\in \Aut(\C/\Q)$.
\endproof

\hfill

\proposition
Let $V_\Q$ be a rational, polarized Hodge structure, $\rho$
the corresponding $U(1)$-action, and $\MT$ its 
special Mumford-Tate group.
Then $\MT$ is a Zariski closure of a group generated by
$\sigma(\rho)$, for all $\sigma\in \Aut(\C/\Q)$.

\hfill

\proof
Since $\MT$ is rational, it is preserved by $\Aut(\C/\Q)$.
The algebraic closure of the group generated by
$\sigma(\rho)$ coincides with $\MT$, because it is 
the smallest group which is rational, Zariski closed 
and contains $\rho$.
\endproof


\section{Transcendental Hodge algebra}


\subsection{Transcendental Hodge lattice}

\definition\label{_transce_latti_Definition_}
Let $V_\Q$ be a rational, polarized Hodge structure of
weight $w$, $V_C= \bigoplus_{ p+q=w \atop p,q \geq 0} V^{p, q}$,
$V^{tr}\subset V_\C$ a minimal Hodge substructure
containing $V^{w,0}$. We call $V^{tr}$ {\bf the transcendental
  Hodge sublattice} of $V_\C$, or {\bf transcendental Hodge lattice}.\footnote%
{Let $V_\Q$ be a rational, polarized Hodge structure of
weight $w$, 
and $V'\subset V_\C$ a minimal rational subspace
containing $V^{w,0}$. In the earlier versions of this
paper I claimed that $V'$ is a Hodge substructure, but this
is false. The space $V'$ is generated by
vectors $v_\sigma\in \sigma(V^{w,0})$ which are weight $w$ eigenvectors of 
operators $\sigma(\rho)$, where $\rho$ denotes the $U(1)$-action,
and $\sigma\in \Aut(\C/\Q)$. However, the set of such operators
is not closed under multiplication, hence $V'$ is not
preserved by the Mumford-Tate group which is generated by those
operators. I am grateful to
Pierre Deligne who noticed this error after this
paper was already published.}

\hfill

\theorem
Transcendental Hodge lattice is a birational invariant.

\hfill

\proof Let $\phi:\; X \arrow Y$ be a birational morphism of projective
varieties. Then $\phi^*:\; H^d(Y)\arrow H^d(X)$ induces
isomorphism on $H^{d,0}$. Therefore, it is injective on 
$H^d_{tr}(Y)$. Indeed, its kernel is a Hodge substructure of
$H^d_{tr}(Y)$ not intersecting  $H^{d,0}$, which is impossible.
Applying the same argument to the dual map, we obtain
that $\phi^*$ is also surjective on $H^d_{tr}(Y)$.
\endproof

\subsection{Transcendental Hodge algebra: the definition}
\label{_tra_alge_Subsection_}

\proposition
Let $M$ be a projective K\"ahler manifold, 
\[ H^*_{tr}(M):=\oplus_d H^d_{tr}(M)\] the direct sum of all
transcendental Hodge lattices, and $H^*_{tr}(M)^\bot$
its orthogonal complement with respect to the polarization
form \[ h(x,y):= \int_M x\wedge y \wedge \omega^{\dim_\C M -\tilde x -\tilde y},
\]
where $\tilde z=i$ for any $z\in \Lambda^i(M)$.
Then $H^*_{tr}(M)^\bot$ is an ideal in the cohomology algebra.

\hfill

\proof
The space $V^*:=H^*_{tr}(M)^\bot$ is a maximal Hodge structure contained in 
\[ A^* :=\bigoplus\limits_{ p+q=w \atop p,q > 0} V^{p, q}.
\] The space $A^*$ is clearly
an ideal. For any two Hodge substructures $X, Y\subset H^*(M)$,
the product $X\cdot Y$ also rational and $U(1)$-invariant, hence
it is also a Hodge substructure. However, $A^*$ is an ideal,
hence $X \cdot V^*$  is a Hodge structure contained in $A^*$.
Therefore, $H^*(M) \cdot V^*$ is contained in $V^*$.
\endproof

\hfill

\definition The quotient algebra $H^*(M)/H^*_{tr}(M)^\bot= H^*_{tr}(M)$
is called {\bf the transcendental Hodge algebra} of $M$.

\hfill

\proposition 
Transcendental Hodge algebra is a birational invariant.

\proof Same as for transcendental Hodge lattices.
\endproof


\section{Zarhin's results about Hodge structures of K3
  type}


\subsection{Number fields and Hodge structures of K3 type}

In this subsection we give a survey of Zarhin's
 results on simple Hodge structures of K3 type given
 in a beautiful and very enjoyable paper \cite{_Zarhin:Hodge_K3_}.
 We use these results in Subsection \ref{_transcende_alge_Subsection_}.

\hfill

\definition
A polarized, rational Hodge structure 
$V_\C= \bigoplus_{ p+q=2 \atop p,q \geq 0} V^{p, q}$
of weight 2 with $\dim V^{2,0}=1$ is called
{\bf a Hodge structure of K3 type}.

\hfill

%
%
%

\remark Let $M$ be a projective K3 surface, and
$V_\Q$ its transcendental Hodge lattice. Then $V_\Q$ is
simple and of K3 type.

\hfill

\proposition (Zarhin) 
Let $V_\Q$ be a simple Hodge structure of K3 sype, and 
$E=\End(V_\Q)$ an algebra of its endomorphisms in the
category of Hodge structures. Then $E$ is a number field.

\hfill

\proof By Schur's lemma, $E$ has no zero divisors, hence
it is a division algebra. Since $E\subset \End(V_\Q)$,
it is countable. To prove that it is a number field, it
remains to show that $E$ is commutative. 
However, $E$ acts on a 1-dimensional space $V^{2,0}$. 
This defines a homomorphism from $E$ to $\C$, which
is injective, because $E$ is a division algebra.
\endproof

\hfill

\theorem (Zarhin)
Let $V_\Q$ be a simple Hodge structure of K3 type, and
$E:= \End(V_\Q)$ its endomorphism field. Then $E$ is either
totally real (that is, all its embeddings to $\C$ are real)
or is an imaginary quadratic extension of a totally real field $E_0$.

\hfill

For the convenience of the reader, we sketch Zarhin's proof here.

\hfill

{\bf Proof. Step 1:}  Let $a\in E$ be an endomorphism, and $a^*$ 
its conjugate with respect to the polarization $h$. Since the polarization
is rational and $U(1)$-invariant, the map $a^*$ also preserves
the Hodge decomposition. Then $a^*\in E$. Denote the generator of
$V^{2,0}$ by $\Omega$. Then 
$h(a(\Omega), a(\bar\Omega)) = h(a^*a(\Omega), \bar\Omega)>0$,
hence $\frac{a^*a(\Omega)}{\Omega}$ is a positive real number.
Therefore, the embedding $E\hookrightarrow \C$ induced by 
$a \arrow \frac{a(\Omega)}{\Omega}$ maps $a^*a$ to a 
positive real number.

\hfill

{\bf Step 2:} Since the group $\Aut(\C/\Q)$ acts on $V$
by automorphisms, it preserves $E$, hence $[E:\Q]$ is a Galois extension.
This means that $\Aut(\C/\Q)$ acts transitively on embeddings from $E$ to $\C$,
preserving the map $a\arrow a^*$.
Therefore, all embeddings $E\hookrightarrow \C$ map $a^*a$ to a 
positive real number.

\hfill

{\bf Step 3:} The map $a \arrow a^*$ is either a non-trivial
involution or identity. In the second case, all embeddings
$E\hookrightarrow \C$ map $a^2$ to a positive real number,
and $E$ is totally real. In the second case, $\tau(a)=a^*$
is an involution, hence its fixed set $E^\tau=: E_0$
is a degree 2 subfield of $E$, with $\tau$ the generator of
the Galois group $\Gal(E/E_0)$.

\hfill

{\bf Step 4:} Let $r$ be the root of the
corresponding quadratic equation $x^2-u=0$. Then $\tau(r)=-r$,
which gives $-r^2 = r\tau(r)>0$, and $-r^2$ is positive and
real for all embeddings $E\hookrightarrow \C$. Therefore,
$[E:E_0]$ is an imaginary quadratic extension.
\endproof

\hfill

\theorem (Zarhin)
Let $V_\Q$ be a Hodge structure of K3 type, and $E:= \End(V_\Q)$
the corresponding number field. Denote by $SO_E(V)$
the group of $E$-linear isometries of $V$ for $[E:\Q]$
totally real, and by $U_E(V)$ the group of $E_0$-linear
isometries of $V$ for $E$ an imaginary quadratic extension
of a totally real field $E_0$. Then the special Mumford-Tate group $\MT$
is $SO_E(V)$ in the first case, and $U_E(V)$ in the second.

\hfill

\proof The group $\MT$ is obtained as the Weil
restriction of scalars from $E$ to $\Q$
for   $SO_E(V)$ and from $E_0$ to $\Q$ for $U_E(V)$.
See \cite{_Zarhin:Hodge_K3_} 
for original proof, \cite{_Zarhin:Hodge_K3_LNM_}
for an alternative proof using B. Kostant's theorem
and \cite{_Zarhin:small_number_}  for generalizations.
\endproof

\hfill

\definition
An imaginary quadratic extension of the totally real field is called
{\bf CM field} (``complex multiplication''), after 
\cite{_Shimura_Taniyama:CM_}. We shall call
the two cases considered in Zarhin's theorem
``the real case'' (for $SO_E(V)$)
and ``CM case'' for $U_E(V)$.

\subsection{Special Mumford-Tate group for Hodge structures of K3 type}
\label{_MT_Zarhin_Subsection_}

Let $V_\C=V_\Q \otimes_\Q \C$ and 
$V_\C=V_\R \otimes_\Q \R$. We would be interested in the
real and complex Lie groups $SO_E(V_\C)$  and $U_E(V_\C)$,
obtained as a complexification of $MT(V)$.

As a Lie group, $SO_E(V_\C)$ is isomorphic to $SO(\C^n)$, where $n$
is dimension of $V_\Q$ over $E$. To see this, consider the tensor product
\[E\otimes_\Q \C=\underbrace{\C\oplus\C\oplus...\oplus\C}_{\text{$k$ times}},\] 
where $k$ is degree of $E$. The $k$ factors of
$E\otimes_\Q \C$ correspond to all embeddings of $E$ to $\C$.
Let $\sigma_1, ..., \sigma_k$ be all such factors, $\sigma_i:\; E\otimes_\Q \C \arrow \C$,
and denote by $V_{\sigma_i}$ the subset of $V= V_\Q\otimes_\Q \C$
corresponding to $\sigma_i$, 
\[
V_{\sigma_i}= \{v\in V_\Q\otimes_\Q \C\ \ \forall a\in E, \ \ a(x)=\sigma_i(x)\}.
\]
Clearly, $V= \bigoplus_i V_{\sigma_i}$, and 
$SO_E(V)$ embeds to each of $SO(V_{\sigma_i})$ tautologically.
Since $V$ is generated as a $E$-module by each of $V_{\sigma_i}$,
this implies that $SO_E(V)=SO(V_{\sigma_i})$, for each $i$.

The same argument is used to show that
in CM-case, $U_E(V_\R)$ is isomorphic to $U(W)$,
for a complex vector space $W$ obtained as follows. Since $E$
is totally imaginary, 
$E\otimes_\Q \R= \underbrace{\C\oplus\C\oplus...\oplus\C}_{\text{$k/2$ times}}$.
Then $V_\R$ is a sum of $k$ copies of the same complex vector space $W$ 
with the Galois group acting transitively, and $U_E(V_\R)=U(W)$.

The space $W\otimes_\R \C$ can be identified with $W\oplus \bar W$, and
the Hermitian form $h$ identifies $\bar W$ with $W^*$, which are
both isotropic subspaces with respect to $h$.
This implies that $U_E(V_\C)$ is a group of isometries of $W\oplus W^*$
preserving this direct sum decomposition, and is identified with
$GL(W)$.


\section{Transcendental Hodge algebra for hy\-per\-k\"ah\-ler manifolds}


\subsection{Irreducible representations of $SO(V)$}

Before we start describing the transcendental Hodge algebra
for a hyperk\"ahler manifold, we have to define a certain
$SO(V)$-invariant quotient algebra of the symmetric algebra 
$\Sym^*(V)$ of a vector space equipped with a non-degenerate
scalar product.

\hfill

\claim\label{_Sym_+-defi_Claim_}
Let $V$ be a vector space over a field $k$, $\Char k=0$
equipped with a non-degenerate symmetric
product $h$, and $b \in \Sym^2(V)$ an $SO(V)$-invariant bivector dual
to $h$. Let $\Sym_+^*(V)$\footnote{We denote $\Sym_+^*(V)$ by
$\Sym_+^*(V)_k$ when $k$-dependence is necessary}
be the quotient of $\Sym^*(V)$ by the
ideal generated by $b$. 
Then $\Sym_+^i(V)$ is irreducible as a representation
of $SO(V)$.

\hfill 

\proof
It is well known (see \cite[Theorem V.5.7.F]{_Weyl:invariants_}) that $\Sym^i(V)$ is decomposed
to a direct sum of irreducible representations of $SO(V)$ as follows:
\[ \Sym^i(V)= \Sym^i_i(V)\oplus \Sym^i_{i-2}(V) \oplus ... \oplus 
\Sym^i_{i-2\lfloor i/2\rfloor }(V),\] where 
$\Sym^i_{i-k}(V)\cong \Sym^{i-k}_{i-k}(V)$ and the subrepresentation
$\Sym^i_{i-2}(V) \oplus ... \oplus \Sym^i_{i-2\lfloor i/2\rfloor }(V)$
is the image of $\Sym^{i-2}(V)$ under the map
$\alpha \arrow \alpha \cdot b$. Then $\Sym_+^i(V)=\Sym_i^i(V)$,
hence irreducible.  \endproof

\hfill

\remark
Since the category of representations of $SO(V)$ is
semisimple, any irreducible quotient of $\Sym^i(V)$ can be 
unambiguously realized as a subrepresentation of $\Sym^i(V)$.
Further on, we shall always consider  $\Sym_+^i(V)$
as a subspace in $\Sym^i(V)$.

\subsection{Transcendental Hodge algebra for hyperk\"ahler manifolds}
\label{_transcende_alge_Subsection_}

\theorem\label{_hyperka_cohomo_1995_Theorem_}
(\cite{_Verbitsky:cohomo_,_Verbitsky:coho_announce_})
Let $M$ be a maximal holonomy hyperk\"ahler manifold, $\dim_\C M=2n$,
and $H^{*}_{(2)}(M, \Q)$ subalgebra in cohomology generated by $H^2(M,\Q)$.
Then $H^{2i}_{(2)}(M,\Q)= \Sym^i H^2(M,\Q)$ for all $i \leq n$.
\endproof

\hfill

\remark
Let $\Omega\in H^2(M, \C)$ be the cohomology class of the
holomorphic symplectic form. By Bogomolov's theorem
(\cite{_Bogomolov:decompo_}), for any hyperk\"ahler
manifold of maximal holonomy the space $H^{*, 0}(M)$ is spanned
(as a vector space) by the powers $\Omega, \Omega^2, ...,
\Omega^n$. In particular, the transcendental Hodge
lattice $H^i_{tr}(M)$ is non-zero only for the even $i$.

\hfill

The main result of this paper:

\hfill

\theorem\label{_transcende_Theorem_}
Let $M$ be a projective
maximal holonomy hyperk\"ahler manifold, $\dim_\C M=2n$, 
$H^{2*}_{tr}(M)$ its transcendental Hodge algebra, 
and $E=\End(V)$ the number field of endomorphisms of its transcendental Hodge lattice
$V= H^2_{tr}(M)$.
Let $\Sym^*(V)_E$ denote the $E$-linear symmetric power.  Then
\[\bigoplus_{i=0}^n  H^{2i}_{tr}(M)= \bigoplus_{i=0}^n
\Sym^i(V)_E\] for $E$ 
a CM-field, and \[ \bigoplus_{i=0}^n  H^{2i}_{tr}(M)= \bigoplus_{i=0}^n \Sym^i_+(V)_E\]
for $E$ totally real.

\hfill

\proof
Using the polarization form on $V_\Q$, we can identify
$\Sym^*_\Q(V_\Q)$ with the space of polylinear symmetric
forms on $V_\Q$. Since $\Omega$ is a common eigenvector
of all elements of $E\subset \End_\Q(V_\Q)$, 
the complexification of $\Sym^*_E(V_\Q)$ contains
$\Omega^i$, and therefore $H^{2i}_{tr}(M)$ belongs
to $\Sym^*_E(V_\Q)$.

By definition, $H^{2i}_{tr}(M)$ is the smallest Hodge substructure
(that is, the special Mumford-Tate subrepresentation)
of $H^{2i}(M,\C)$ containing $\Omega^i$. By \ref{_hyperka_cohomo_1995_Theorem_},
one could replace $H^{2i}(M,\C)$ with $\Sym^i(V_\C)$.
This means that we need to find the smallest $U_E(V)$-
and $SO_E(V)$-representations in $\Sym^i(V_\C)$ containing
$\Omega^i$. In CM-case $V$ is the standard representation of $U(V)$, 
the space $\Sym^i(V)$ is irreducible. 

For the totally real
case , Mumford-Tate group is $SO_E(V)$. However,
the bivector $b$ of \ref{_Sym_+-defi_Claim_} is of Hodge
type $(2,2)$, and the representations 
$\Sym^i_{i-k}(V)$ defined in the proof of 
\ref{_Sym_+-defi_Claim_} are orthogonal to $\Omega^i$. 
Therefore, the form $\Omega^i$ 
belongs to  $\Sym^i_+(V)_E$, which is an irreducible
representation of $SO_E(V)$   (\ref{_Sym_+-defi_Claim_}) identified with the
transcendendental Hodge lattice.
\endproof

\hfill

\definition
Let $S$ be a family of deformations of a hyperk\"ahler
manifold, $s\in S$ a point, $M_s$ its fiber, and $T_sS
\arrow H^1(TM_s)$ the corresponding tangent map, defined
by Kodaira and Spencer. The dimension of its image in
general $s\in S$ is called {\bf the essential dimension}
of the deformation family. 
Essential dimension of a holomorphic family $S$ of
deformations is a smallest $d$ such that $S$
can be locally obtained as a pullback of a $d$-dimensional
family.

\hfill

The dimension of the transcendental Hodge lattice of a general
member of a family is determined
by the essential dimension of the deformation family:

\hfill

\theorem\label{_MT_family_Theorem_}
Let $M$ be a projective hyperk\"ahler manifold, which is generic in
a family $S$. Assume that the essential dimension 
of this deformation family is $d$. Then
$\dim_E H^2_{tr}(M)\leq d+2$, where $E=\End(H^2_{tr}(M))$
is the corresponding number field.

\hfill

{\bf Proof:} Using the Torelli theorem (\cite{_V:Torelli_}), we may 
assume that $S$ is a subset of the period space $\Perspace$ of 
the Hodge structures on $H^2(M)$.
Then \ref{_MT_family_Theorem_} 
is implied by \ref{_dimension_family_MT_Proposition_} below. \endproof


\section{Variations of Hodge structures of K3 type and 
the special Mumford-Tate group}


\subsection{Variations of Hodge structures and the special Mumford-Tate group}

\definition
Let $M$ be a complex manifold, and $V$ 
a real vector bundle equipped with a flat connection
over $M$. Assume that at each point $m\in M$, the fiber $V_m$
is equipped by a rational, polarized Hodge structure, that is,
a Hodge decomposition $V_m \otimes \C = \bigoplus V_m^{p,q}$,
smoothly depending on $M$, a
rational structure $V_m=V_m(\Q)\otimes_\Q \R$ and a polarization
$s\in V_m(\Q)^*\otimes V_m(\Q)^*$. 
Suppose that the following conditions are satisfied.
First, the rational lattice
$V_m(\Q)$ and the polarization are preserved by the connection.
Second, the Hodge decomposition $V_m \otimes \C = \bigoplus V_m^{p,q}$
satisfies {\bf the Griffiths transversality condition}:
\[
\nabla(V^{p,q})\subset \bigg(V^{p,q}\otimes \Lambda^1(M)\bigg) \oplus 
\bigg(V^{p+1,q-1}\otimes \Lambda^{0,1}(M)\bigg) \oplus 
\bigg(V^{p-1,q+1}\otimes \Lambda^{1,0}(M)\bigg). 
\]
Then $(V, V_m \otimes \C = \bigoplus V_m^{p,q}, V_m(\Q), s)$
is called {\bf a rational, polarized variation of Hodge structures},
abbreviated to VHS.

\hfill

\example
Let  ${\cal X}\arrow S$ be a holomorphic family of 
compact K\"ahler manifolds over a base $S$, and $V$ the corresponding vector bundle
of cohomology of the fibers ${\cal X}_t, t\in S$. Clearly,
$V$ is equipped with the tautological connection $\nabla$,
called {\bf the Gauss-Manin connection}. The Hodge decomposition 
on $H^*({\cal X}_t)$ depends smoothly on $t\in S$. If, in addition,
all ${\cal X}_t$ admit a polarization (that is, an ample line
bundle) $L_t$ with $c_1(L_t)$ invariant with respect to $\nabla$. Then
the sub-bundle of primitive forms is also $\nabla$-invariant, and the
corresponding Hodge structures glue together to give 
a polarized, rational
variation of Hodge structures (see \cite{_Griffiths:transcendental_} or \cite{_Voisin-Hodge_}).

\hfill

From \ref{_Mumford_Tate_via_invariants_Theorem_}
it is obvious that the special Mumford-Tate group is lower
semicontinuous in smooth families of Hodge structures. 
Indeed, it is the biggest group
which fixes all Hodge vectors in all 
tensor powers of $V$. The sets Hodge vectors are upper
semicontinuous as functions of a base, 
because for each rational vector, the
set of fibers where it is of Hodge type $(p,p)$ is closed.
This set has a special name in the theory of variations
of Hodge structures.

\hfill

\definition
Let $V$ be a variation of Hodge structures over $S$,
and $\phi\in V^{\otimes n}$ any tensor. The set of
all $x\in S$ such that $\phi$ is of Hodge type $(p,p)$
in $x$ is called {\bf the Hodge locus} of $\phi$.

\hfill

The following claim is easily deduced from Griffiths
transversality condition.

\hfill

\claim  (\cite{_Griffiths:transcendental_}, \cite{_Cattani_Kaplan:algebraicity_}
\cite[II, 5.3.1]{_Voisin-Hodge_})\\
Let $V$ be a variation of Hodge structures over a base
$S$, $\phi\in V^{\otimes n}$ any tensor, and $S_\phi$
the corresponding Hodge locus. Then $S_\phi$ is a complex analytic 
subset of $S$.
\endproof

\hfill

\corollary\label{_MT_generic_Corollary_}
Let $V$ be a VHS over a connected base $S$, and $R$ the set of all 
rational vectors $\phi\in V_\Q^{\otimes n}$  such that
the corresponding Hodge locus $S_\phi$ is a proper subset of $S$,
and $S_0\subset S$ the complement $S\backslash \bigcup_{\phi\in R} S_\phi$.
Then for each $t\in S_0$ the special Mumford-Tate group $MT(V_t)$ 
is independent from $t$, and contains $MT(V_s)$ for any
other $s\in S$. 

\hfill

{\bf Proof:} Follows immediately from
\ref{_Mumford_Tate_via_invariants_Theorem_}. \endproof

\hfill

\definition
In assumptions of \ref{_MT_generic_Corollary_},
a point $s\in S_0$ is called a {\bf Mumford-Tate generic} point of $S$.

\subsection{Variations of Hodge structures of K3 type}

\remark\label{_periods_K3_type_Remark_}
Consider the space $\Perspace$ of all K3-type Hodge structures
on $(V_\Q, q)$, where $V_\Q$ is a rational vector space and
$q$ a rational bilinear symmetric form on $V_\Q$ of signature $(2,n)$.
Let $V_\C:= V_\Q \otimes_\Q \C$. Then 
\[
\Perspace = \{ l\in {\Bbb P} V_\C\ \ |\ \ q(l,l)=0, q(l, \bar l) >0\}.
\]
Indeed, for each $l$ in this set, we can define a decomposition
$V^{2,0}=l$, $V^{0,2}=\bar l$, $V^{1,1}=\langle V^{2,0}+ V^{0,2}\rangle^\bot$.
The following claim is well known.

\hfill

\claim
The trivial local system $(V_\Q\otimes_\Q \R, q)$ over $\Perspace$
with the Hodge decomposition defined in \ref{_periods_K3_type_Remark_}
is a VHS.

\hfill

{\bf Proof:} The only non-trivial part is Griffiths transversality:
we need to show that  
\begin{equation}\label{_Griffiths_tra_for_K3_type_Equation_}
\nabla(V^{2,0})\subset V^{2,0}\otimes \Lambda^1(M) 
\oplus V^{1,1}\otimes \Lambda^{1,0}(M). 
\end{equation}
Choose any section $\xi$ of the line bundle $V^{2,0}$
on $\Perspace$. Differentiating $q(\xi, \xi)=0$, 
we obtain $q(\nabla(\xi), \xi)=0$, 
which implies \eqref{_Griffiths_tra_for_K3_type_Equation_}
immediately. \endproof

\hfill

\proposition\label{_dimension_family_MT_Proposition_}
Let $V$ be a VHS of K3 type over a complex manifold $S\subset \Perspace$,
and $G$ the special Mumford-Tate group of a generic point $s\in S$. 
Then $G\cong SO_E(V_\Q)$ or $G\cong U_E(V_\Q)$ (Subsection
\ref{_MT_Zarhin_Subsection_}), and $\dim S\leq \dim_E V_\Q-2$.

\hfill

{\bf Proof:} Let $W\subset T^{\otimes} V$ be the set of all
rational tensors $\phi \in T^{\otimes} V$ which are of type
$(p,p)$ in $s$. By \ref{_Mumford_Tate_via_invariants_Theorem_},
$W$ is the space of $G$-invariants: $W= (T^{\otimes} V)^G$.
Replacing $S$ by a bigger complex analytic subvariety of $\Perspace$, we may assume
that $S$ is an open subset of the intersection $S_1$ of all Hodge loci for all
$w\in W$; indeed, this intersection is complex analytic,
contains $S$, and has the same generic special Mumford-Tate group.

The space $S_1$ is preserved by $G$ acting on ${\Bbb P}V_\C$;
using \ref{_Mumford_Tate_via_invariants_Theorem_}
again, we identify $S_1$ with 
an orbit of $G$. For $G=SO_E(V)$, this orbit
has complex dimension $\dim_E V_\Q-2$, because it contains
an open subset of a projectivization of an appropriate quadric,
and a quadric in ${\Bbb P} A$ has dimension $\dim A-2$.

For $G=U_E(V)$, we have $U_E(V_\C)= GL(W)$ (Subsection 
\ref{_MT_Zarhin_Subsection_}) acting on $V=W\oplus W^*$.
Clearly, $S$ contains a point on a quadric 
\[ Q=\{l\in {\Bbb P}V \ \ |\ \ l=\C \cdot (x_1, x_2)\ \ |\ \ \langle x_1, x_2\rangle =0\}.\]
The whole quadric is an orbit of $GL(W)=U_E(V_\C)$, hence
it contains $S$, and its dimension is again $\dim_E V_\Q-2$.
\endproof


\section{$k$-symplectic structures and symplectic embeddings}


\subsection{Non-degeneracy of the map $x\arrow x^n$ in $H^*_{tr}(M)$.}

Applications of \ref{_transcende_Theorem_} are based on the 
following observation.

\hfill

\theorem
Let $M$ be a projective 
maximal holonomy hyperk\"ahler manifold, $\dim_\C M=2n$ and
$x\in H^2_{tr}(M)$ a non-zero vector. 
Then $x^n\neq 0$ in $H^n_{tr}(M)$.

\hfill

\proof
When $E$ is a CM-field and $H^*_{tr}(M)=\Sym^*(V)$, this is trivial.
When $E$ is totally real and $H^*_{tr}(M)=\Sym^*_+(V)$, this is proven as follows.
The action of $SO(V)$ on the projectivization ${\Bbb P}V$
has only two orbits: the quadric $Q:= \{x \ \ |\ \ h(x,x)=0\}$
and the rest. The map $P(x) = x^n$ is non-zero, hence it is
non-zero on the open orbit. To check that it is non-zero
on the quadric, notice that $\Omega$ belongs to $Q$, and
$\Omega^n$ is non-zero in $\Sym^n_+(V)$.
\endproof

\hfill

\corollary\label{_symple_torus_non-dege_Corollary_}
Let $M$ be a  projective 
maximal holonomy hyperk\"ahler manifold, 
$N \hookrightarrow M$ a holomorphic symplectic embedding, $\dim_\C N = 2n$,
and $x\in H^2(N)$ a restriction of a non-zero class 
$x\in H^2_{tr}(M)$. Then $x^n \neq 0$.

\hfill

\proof Since the embedding $N\hookrightarrow M$ 
is holomorphically symplectic, the restriction map
$H^{2n}_{tr}(M)\arrow H^{2n}_{tr}(N)$ is non-zero,
hence injective. On the other hand, $x^n \neq 0$
in $H^{2n}_{tr}(M)$, as shown above.
\endproof

\subsection{$k$-symplectic structures: definition and applications.}

\definition 
Let $V$ be a $4n$-dimensional vector space, $W$ a $k$-dimensional 
vector space, and $\Psi:\; W \arrow \Lambda^2(V)$ a linear map. Assume that
$\Psi(\omega)$ is a symplectic form for general $\omega\in W$, and has rank 
$\frac 1 2 \dim W$ for $\omega$ in a non-degenerate quadric $Q\subset W$. 
Then $\Psi$ is called {\bf $k$-symplectic structure on $V$}.

\hfill

\remark It is not hard to see that the set of degenerate $\alpha\in W$ 
is always a quadric, if the rank of degenerate $\alpha$ is $\frac 1 2 \dim W$
(see \cite{_SV:k-symplectic_}), but this quadric can be degenerate.

\hfill

\theorem \label{_k-symple_Clifford_Theorem_}
(\cite{_SV:k-symplectic_})
Let $V$ be a $k$-symplectic space.
Then $V$ is a Clifford module over a Clifford algebra $\Cl(W_0)$,
with $\dim W_0= \dim W-1=k$, and 
$\dim V$ is divisible by $2^{\left\lfloor(k-1)/2\right\rfloor}$.
\endproof

\hfill

\corollary \label{_k_symple_in_families_Corollary_}
Let $M$ be a projective hyperk\"ahler manifold, 
generic in a deformation family of dimension $d$,
and $N \rightarrow M$ a compact torus which is holomorphically symplectically
immersed to $M$.
Then $\dim N$ is divisible by $2^{\left\lfloor \frac{d+1}2\right\rfloor}$.

\hfill

\proof
By \ref{_symple_torus_non-dege_Corollary_}, $H^1(T)$ is $k$-symplectic, where 
$k = \dim_E H^2_{tr}(M)$, and $E=\End(H^2_{tr}(M)$ the number
field associated with $M$. Indeed, fix an embedding $E\hookrightarrow \C$,
ad let $A:=H^2_{tr}(M)\otimes \C$. This space is equipped with 
a natural action of the complexification of the special Mumford-Tate
group $G_\C$, and $A$ has only two non-zero orbits
with respect to $G_\C$-action (Subsection \ref{_MT_Zarhin_Subsection_}).
This implies that non-zero 2-forms on $H^1(T)$ induced from $A$
can possibly be non-degenerate or have rank $r< \dim H^1(T)$ for
a fixed $r>0$. This implies that for any two symplectic forms
 $\omega, \omega'\in A$, the operator 
$\omega^{-1} \circ \omega':\; H^1(T,\C)\arrow H^1(T,\C)$
has at most two different eigenvalues, and the corresponding
eigenspaces have the same dimension: $r= \frac {\dim H^1(T)}{2}$.
The space of degenerate 2-forms $\eta \in A$ is a quadric, and
they all have rank either 0 or $r$, hence $A$ defines a $k$-symplectic
structure on $H^1(T, \C)$.

Then $k\geq d+2$, 
as follows from \ref{_MT_family_Theorem_}. By 
\ref{_k-symple_Clifford_Theorem_} it is a Clifford module
over $\Cl_E(d+1)$, and 
$\dim H^1(T)$ is divisible by $2^{\left\lfloor {d+1}/2\right\rfloor}$.
\endproof

\hfill

This means, in particular, that 2-dimensional symplectic
tori exist only in families of essential dimension $\leq 3$
and 4-dimensional symplectic
tori exist only in families of essential dimension $\leq 5$.
Also, for a general projective hyperk\"ahler manifold with
$b_2=23$  and $\Pic=\Z$, any symplectic subtorus has
dimension divisible by $2^{10}=1024$, because the
corresponding deformation space has dimension $b_2(M)-3=20$.

\hfill

As another application, we generalize 
\cite[Corollary 1.17]{_SV:k-symplectic_}.

\hfill

\theorem
Let $Z\subset M$ be a symplectic subvariety in a
hyperk\"ahler manifold. Then $H^2_{tr}(M)$ maps
injectively to $H^2_{tr}(Z)$. In particular, if 
the pair $(M,Z)$ can be deformed in a family $S$
such that the essential dimension
of the corresponding deformation of $M$ is $d$,
then $\dim H^2_{tr}(Z)\geq (d+2) e$, where
$e=\dim_\Q E$ and $E=\End_\Q (H^2_{tr}(M))$
is the corresponding number field.

\hfill

{\bf Proof:} The first statement is clear because
$H^2_{tr}(M)$ is irreducible, and mapped to $H^2_{tr}(Z)$
non-trivially because $Z$ is holomorphically symplectic. 

To prove the second statement, notice that 
dimension of $H^2_{tr}(M)$ over $E$ is at least $d+2$
(\ref{_MT_family_Theorem_}), hence the dimension of
$H^2_{tr}(M)$ over $\Q$ is at least $(d+2) e$,
and this space is embedded to $H^2(Z,\Q)$.
\endproof

\hfill

{\bf Acknowledgements:}
I am grateful to Yuri Zarhin for generously sharing his
insight in the transcendental Hodge lattices and to
Ljudmila Kamenova for many interesting discussions.
Million thanks to Pierre Deligne who found an error
in the published version of this paper (see the footnote
to \ref{_transce_latti_Definition_}).

\hfill

{\small
\noindent {\sc Misha Verbitsky\\
{\sc Universit\'e Libre de Bruxelles, CP 218,\\
Bd du Triomphe, 1050 Brussels, Belgium\\{
\tt  mverbits@ulb.ac.be}, \\
also: \\
{\sc Laboratory of Algebraic Geometry,\\
National Research University Higher School of Economics,\\
Department of Mathematics, 6 Usacheva street, Moscow, Russia.}\\
}}}

\end{document}